\documentclass{elsarticle}
\usepackage{amsmath}
\usepackage{amssymb}
\usepackage{graphics}
\usepackage{epstopdf}
\usepackage{float}
\newtheorem{theorem}{Theorem}[section]

\begin{document}

\begin{frontmatter}

\title{Analysis of Impact Chattering}

\author[akh]{M. ~Akhmet\corref{cor1}}
\ead{marat@metu.edu.tr}

\author[akh]{S. ~\c{C}a\u{g}}
\ead{sabahattincag@gmail.com}

\cortext[cor1]{Corresponding author}

\address[akh]{Department of Mathematics, Middle East Technical University, Ankara, Turkey.}

\begin{abstract}
In this paper, mechanical models  with    Newton's Law  of impacts are studied. One of the most interesting properties in some of these models is chattering. This phenomenon is understood as the appearance of an infinite number of impacts occurring in a finite time. Conclusion  on  the presence of chattering  is   made  exclusively    by  examination of  the   right hand   side  of impact  models for the first time. Criteria for the sets of initial data which always lead to chattering are established. The Moon-Holmes model is subject to regular impact perturbations for the chattering generation.
Using the chattering solutions, continuous chattering is generated.  To  depress  the chattering, Pyragas control is applied.  Illustrative examples are provided to demonstrate the impact chattering.
\end{abstract}
 
\begin{keyword}
Impact  mechanism, The Moon-Holmes model, Chattering, Bouncing bead, Pyragas control.
\end{keyword}

\end{frontmatter}

The implementation of sliding mode control is often irritated by high frequency oscillations known as ``chattering'' in system outputs issued by dynamics from actuators and sensors ignored in system modeling \cite{Lee}. In study \cite{ibrahim}, chattering is considered as a special type of oscillation characterized by very small amplitudes that are decreasing with time. In impacting systems, it is understood as an infinite number of discontinuities moments occurring in a finite time period, for instance, a ball bouncing to rest on a horizontal surface \cite{Giusepponi}. It is asserted in \cite{Giusepponi} that chattering resembles with the inelastic collapse. The balls dissipate their energy through an infinite number of collisions in a finite time interval.  Budd and Dux \cite{budd} showed that chattering can occur for a periodically forced, single degree of freedom impact oscillator with a restitution law. They demonstrated that chattering can form part of a periodic motion, and this  relates to certain types of chaotic behavior. However, they studied through an example. Using the solution, they proved the existence of chattering for a linear system.

Nordmark and Piiroinen\cite{nordmark} considered simulation problems for chattering as well as analysis of stability of the limit cycle, which is chattering by solving the first variational equations. Moreover, they used the mappings, which are constructed  with the help of a solution, in simulation schemes. Similar to the one  in paper \cite{budd}, it was shown that the existence of chattering for a linear system.  Nonetheless, in both papers \cite{budd,nordmark}, they do not consider the conditions which guarantee the appearance of chattering. In this study, we consider the chattering as a motion with infinite number of discontinuities  in a finite time. This is the first time that  sufficient conditions are provided for the chattering based on  properties not on maps derived with the help of solutions, but, on conditions for the right-hand side of impulsive systems. Our models essentially are nonlinear (see, for example, Example 1). Since this is the first result in this direction, the models under consideration are respectively simple. Nevertheless, this is a class of mechanical models which can be significantly enlarged in the future investigations by consideration of large ensembles of impact oscillators and weakening conditions of the present paper. We consider models with vibrating surface of impacts as well as analyzed problems of Pyragas controllability and existence  of continuous chattering for a model connected unilaterally to a system with an impact chattering. An interesting problem of the regular perturbation of a system with chattering is discussed.

A particular feature of system with impacts  is the existence of the chattering. We have two different types of it, namely complete and incomplete chattering \cite{budd,nordmark}. Complete chattering is the phenomenon wherein a system an infinite number of discontinuities in a finite time occurs, where the velocity tends to zero uniformly. Incomplete chattering bears on a sequence of the impacts that initially has the same behavior as complete chattering, but it ends after a  large but finite number of impacts \cite{nordmark}. In section 2, we will discuss the transient chattering for systems with small parameter considering the transformation of the incomplete chattering to the complete one when the parameter diminishes to zero.

It was first  found by Arnold \cite{Arnold}  that the significant characteristic property of chatter vibration is that it is not generated by external periodic forces, but rather it is generated in the dynamic process itself. Therefore, it is important to emphasize that the systems under investigation in this paper are autonomous.

Consider the problem of  impact interaction of a body falling in the uniform gravity force field with a fixed horizontal base. After colliding with the base the body bounces back with the velocity whose norm is equal to the norm of the pre-impact velocity multiplied by $r,$ where $r$ is the restitution coefficient, $0<r<1.$ Then, after some time interval the body will fall on the base again and the norm of its velocity will be equal to the norm of bouncing velocity in the previous collision multiplied by $r.$ The process cannot end in a finite number of collisions.
Thus, the considered phenomenon consists in following: after the initial collision a series of repeated collisions of attenuated to zero, which ends in a finite time with establishing a long contact between interacted bodies. Arising this contact results in decreasing number of degrees of freedom of the system by a unit or more. So, it is reasonable to call this phenomenon the impact chattering.

It is shown by investigations and observations  that the impact chattering meets in operating almost every mechanism and machine of impact-oscillating type \cite{nagaev}.  Various problems of impact chattering are far from trivial, and their solutions cannot be obtained in closed form for rather general case. As for the use of approximate analytical and numerous methods, it is simplified essentially if one proceeds from the conception about infinity number of impacts inside a finite time range. For example, the  existence of impact chattering was investigated in \cite{nagaev, Kob}.  They simply consider  the free falling of a bead on an immobile base  and on a vibrating table with constant velocity.  In this paper, we consider a more general system and prove the existence of impact chattering.

The chattering phenomena are unwanted in engineering since it is an appearance of multi-strikes in a short period of time. It is not desirable in models of mechanics since it appears as infinite discontinuities in a finite time which make theoretical analysis difficult. We have a research plan to consider theoretical and mathematical complexities connected to chattering and we approach the problem from one of the two possible  points of view. The first one is when mechanical models changed  such that the theoretical chattering disappears \cite{akhmetayse}. The other point of view, which is considered in this paper, is that we approximate a model with infinite moments of discontinuities with those having a finite number of impacts.

This article is organized as follows. The impact model is stated in the first section. In this model each collision is assumed instantaneous, and it comes to rest after an infinite
number of impulse moments in a finite time. The existence of chattering is proved. Asymptotic approximation of solutions with chattering are discussed in section 2. Then, we show that the chattering occurs for a bead bouncing on a sinusoidally vibrating table in section 3. The modified Moon-Holmes model with a small perturbation is discussed in section 4. Using the continuous dependence on parameters and initial value for the impulsive differential equations with non-fixed moments, it is shown  that the solution of the modified Moon-Holmes model is chattering. Following that, the appearance of  continuous chattering by perturbation method is demonstrated in section 5. Finally, by Pyragas control method the chattering solution is controlled to be  periodic.

\section{Existence of Chattering}
An impacting system  admits  a chattering  if there    is a   solution   with  infinite  impulse moments in a finite time.  Moreover,  we   will   say  that   a perturbed   system   admits a transient chattering, if a number  of impacts increases to  infinity  on a  fixed interval as the  small  parameter tends to  zero.

A mechanism with a rigid flat surface  of impacts and the constant coefficient of restitution $r,$ $0<r<1,$ can be modeled by the
following impulsive system
\begin{equation}\label{mainc}
\begin{split}
	\ddot{x}&=f(x,\dot{x}), \\
	\Delta \dot{x}|_{x=\varphi}&=-(1+r)\dot{x}, 
\end{split}
\end{equation}
where $x(t)$ is the coordinate of the bead which is over the impact surface  $x=\varphi, $ $\dot{x}(t)$ is its velocity, $f(u,v)$ is a continuous function on the domain $H=\{0< \varphi \leq u \leq h, |v|\leq \bar{h} \}$ for fixed positive numbers $h, \bar{h},$ and it satisfies the local Lipschitz  condition in its variables on $H$. The equality $\Delta \dot{x}(\theta) = \dot{x}(\theta+)-\dot{x}(\theta-)$ denotes the jump operator in which $t=\theta$ is the time when the bead reaches the rigid obstacle, $\dot{x}(\theta-)$ is the pre-impact velocity and $\dot{x}(\theta+)$ is the post-impact velocity. 

In system \eqref{mainc}, we need the following conditions.
\begin{itemize}
	\item[(C1)] There is a positive number $m$ such that $f(u,v)<-m$ for all $(u,v)\in H$,
	\item[(C2)] $f(u,v)=f(u,-v)$ for all $(u,v)\in H.$
\end{itemize}
 Conditions on function $f(u,v)$ and compactness of domain $H$ imply that there exists a positive number $M$ such that $f(u,v)\geq -M$ for all $(u,v)\in H.$ 
\begin{theorem}\label{thmc}
If conditions (C1), (C2) are satisfied and the following inequality

\begin{equation}\label{conditionhbar}
M\sqrt{\frac{2(h-\varphi)}{m}}<\bar{h}
\end{equation}
is valid, then all solutions with initial value $(x(0),\dot{x}(0))=(x_0,0),$ $\varphi<x_0<h,$  of system \eqref{mainc} are chattering.
\end{theorem}
\textbf{Proof.}
 Consider   an initial value $(x_0,0)\in H,$ $\varphi<x_0<h.$  Denoting $x_1=x, x_2=\dot{x}$    present   the   system \eqref{mainc} as 
\begin{equation}\label{mainc2}
\begin{split}
	\dot{x_1}&=x_2, \\
	\dot{x_2}&=f(x_1,x_2), \\
	\Delta x_2 |_{x_1=\varphi}&=-(1+r)x_2.
\end{split}
\end{equation}
The solution of system \eqref{mainc2} starting at $(x_0,0)$   is 
\begin{subequations}
\begin{equation}\label{solution0a}
x_1(t)=x_0+\int_{0}^t{(t-s)f(x_1(s),x_2(s))}ds,
\end{equation}
\begin{equation}\label{solution0b}
x_2(t)=\int_{0}^t{f(x_1(s),x_2(s))}ds,
\end{equation}
\end{subequations}
while it is continuous.
 By Eq. \eqref{solution0a} and condition $(C1),$   the   coordinate  $x_1(t)$ decreases to $\varphi$ such that there exists a moment $\theta_1$ where  $x_1(\theta_1)=\varphi$ and $x_2(\theta_1)<0.$ Moreover, $x_1(\theta_1+)=\varphi$ and $x_2(\theta_1+)=-rx_2(\theta_1)>0.$ 

Let us show that the solution is continuable to $+\infty$ and it remains in the domain $H.$ First of all, consider   the  interval  $[0,\theta_1].$ From conditions (C1) and (C2), it implies that $x_1(t)\leq x_0 < h, t\in[0,\theta_1].$
Using  \eqref{solution0a} and   inequality  $\varphi<x_0<h$ we get
\begin{equation*}\label{thetaboundforx2}
|h-\varphi|> |\varphi-x_0|= \left| \int_{0}^{\theta_1}{(\theta_1-s)f(x_1(s),x_2(s))}ds \right|\geq \int_{0}^{\theta_1}{(\theta_1-s)m}ds=m\frac{\theta_1^2}{2},
\end{equation*}
which implies that $\theta_1 < \sqrt{\frac{2(h-\varphi)}{m}}.$ 

Consequently, from \eqref{solution0b} and condition \eqref{conditionhbar}
$$|x_2(\theta_1)|=\left| \int_{0}^{\theta_1}{f(x_1(s),x_2(s))}ds \right| \leq M\theta_1<\bar{h}.$$ Thus, we obtain   that $\varphi\leq x_1(t)<h$ and $|x_2(t)|<\bar{h}$ for $t\in [0,\theta_1].$

Applying the   same   arguments as  for  $\theta_1$  one can  show that  there is an  intersection moment $ \theta_2$  such that   $x_1(\theta_2) = \varphi$   and  $x_1(t) > \varphi, t\in (\theta_1,\theta_2).$  In this interval, we have   
\begin{subequations}
\begin{equation}\label{solution3a}
x_1(t)=\varphi+x_2(\theta_1+)(t-\theta_1)+\int_{\theta_1}^t{(t-s)f(x_1(s),x_2(s))}ds,
\end{equation}
\begin{equation}\label{solution3b}
x_2(t)=x_2(\theta_1+)+\int_{\theta_1}^t{f(x_1(s),x_2(s))}ds,
\end{equation}
\end{subequations}
By condition $(C2)$, $x_2(\theta_1+)$ is the maximum value of $|x_2(t)|$ for $t\in (\theta_1,\theta_2].$  Thus, $|x_2(t)| \le r|x_2(\theta_1)|<r\bar{h}<\bar{h}.$ Moreover, from conditions $(C1)$ and $(C2),$ there exists a moment $\xi_1$, $\theta_1<\xi_1<\theta_{2}$, such that $x_2(\xi_1)=0$ and $x_1(\xi_1)$ is the maximum value of $x_1(t)$ on $(\theta_1,\theta_2]$.   Thus, $x_1(t)\leq x_1(\xi_1)<x_0<h,$ 
and the   trajectory   of  $x(t)$ is in $H$ for $t\in [\theta_1,\theta_2].$  Next, recursively, it can be shown that there exists an increasing sequence $\theta_i, i = 1,2,\ldots,$ such that $x_1(\theta_i)=\varphi,i=1,2,\ldots,$   and   the   orbit   of $x(t)$   is  in $H$    for   all $t \ge 0.$   

Now, we   will   show that   the sequence $\theta_i$  converges.  The solution of system \eqref{mainc2}  is defined by 
\begin{subequations}\label{solution1}
\begin{equation}\label{solution1a}
x_1(t)=\varphi+x_2(\theta_i+)(t-\theta_i)+\int_{\theta_i}^t{(t-s)f(x_1(s),x_2(s))}ds,
\end{equation}
\begin{equation}\label{solution1b}
x_2(t)=x_2(\theta_i+)+\int_{\theta_i}^t{f(x_1(s),x_2(s))}ds.
\end{equation}
\end{subequations}
on the interval $(\theta_i,\theta_{i+1}], i =1,2,\ldots.$
 
Using condition (C1), it can be shown that there exists a moment $\xi_i$, $\theta_i<\xi_i<\theta_{i+1}$, such that $x_2(\xi_i)=0$. Also, utilizing condition (C2), we obtain $\xi_i=\frac{\theta_i+\theta_{i+1}}{2}.$   The solution on the interval $(\theta_{i+1},\theta_{i+2}]$ is 
\begin{equation}\label{solution2}
\begin{split}
x_1(t)&=\varphi+rx_2(\theta_i+)(t-\theta_{i+1})+\int_{\theta_{i+1}}^t{(t-s)f(x_1(s),x_2(s))}ds, \\
x_2(t)&=rx_2(\theta_i+)+\int_{\theta_{i+1}}^t{f(x_1(s),x_2(s))}ds.
\end{split}
\end{equation}  
From $x_2(\xi_i)=0$ and $x_2(\xi_{i+1})=0,$ we get 
\begin{eqnarray}
x_2(\theta_i+)=-\int_{\theta_i}^{\xi_i}{f(x_1(s),x_2(s))}ds,\label{x2}\\
rx_2(\theta_i+)=-\int_{\theta_{i+1}}^{\xi_{i+1}}{f(x_1(s),x_2(s))}ds.\label{rx2}
\end{eqnarray}
Let us divide \eqref{rx2} by \eqref{x2} in order to get
\begin{equation*}
\displaystyle r=\frac{\int_{\theta_{i+1}}^{\xi_{i+1}}{f(x_1(s),x_2(s))}ds}{\int_{\theta_i}^{\xi_i}{f(x_1(s),x_2(s))}ds}.
\end{equation*}
Using mean value theorem, we have 
\begin{equation}\label{thetaapprox}
\displaystyle r=\frac{(\xi_{i+1}-\theta_{i+1})f(x_1(s^*),x_2(s^*))}{(\xi_i-\theta_i)f(x_1(s^{**}),x_2(s^{**}))},
\end{equation}
for   some  $s^{**}$ and $s^{*}$   in $(\theta_{i},\theta_{i+1})$  and  $(\theta_{i+1},\theta_{i+2})$  respectively.

Then, 
\begin{equation}
\displaystyle \frac{\theta_{i+2}-\theta_{i+1}}{\theta_{i+1}-\theta_i}=\frac{\xi_{i+1}-\theta_{i+1}}{\xi_{i}-\theta_{i}}<\frac{M_i}{m_i}r, i=1,2,3....,
\label{thetas}
\end{equation}
where  $\displaystyle M_i=\max_{[\theta_{i},\theta_{i+1}]}|f(x_1(t),x_2(t))|$ and $\displaystyle  m_i=\min_{[\theta_{i},\theta_{i+1}]}|f(x_1(t),x_2(t))|.$ 
Since $r<1,$    $\displaystyle \max_{[\theta_{i},\theta_{i+1}]}|x_1(t)|\to \varphi$ and $\displaystyle  \max_{[\theta_{i},\theta_{i+1}]}|x_2(t)| \to 0$ as $i \to \infty.$   Moreover, continuity of $f(u,v)$ implies that $\frac{M_i}{m_i} \to 1$ as $i \to \infty.$ This and \eqref{thetas}  prove   the convergence.  The theorem   is proved.

\begin{figure}[H]
\centering
\includegraphics[width=\columnwidth]{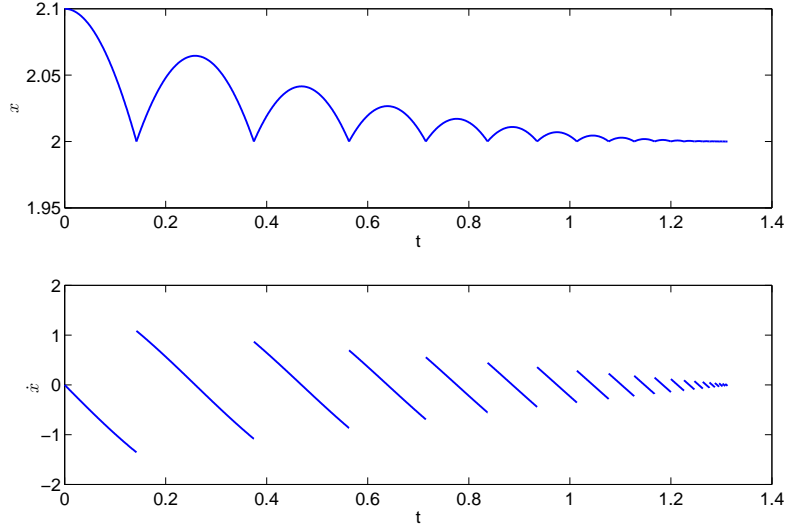}%
\caption{The   graphs of coordinates $x(t)$ and $\dot{x}(t)$ with initials $x(0)=2.1$ and $\dot{x}(0)=0$   of  system \eqref{exc} with $r=0.8.$}  
\label{nlinear}%
\end{figure}

\textbf{Example 1. } Consider the following non-linear   system
\begin{equation}
\begin{split}
\ddot{x}+\cos(\dot{x})+x^3=0, \\
\Delta \dot{x}|_{x=2}=-(1+r)\dot{x},
\end{split}
\label{exc}
\end{equation}
in the domain  $2\leq x \leq 2.5, |\dot{x}|<7.$  We have  $f(x,\dot{x})=-\cos(\dot{x})-x^3\leq -7$,  $f(x,\dot{x})=f(x,-\dot{x})$, $f(x,\dot{x})\geq -16.625$   in the domain.  Condition \eqref{conditionhbar} is true since $16.625\sqrt{1/7}\approx 6.28<7.$ That   is,  we  are  in circumstances   of Theorem  \ref{thmc}   and  if we choose $r=0.8$, $x(0)=2.1,\dot{x}(0)=0,$  the solution of system \eqref{exc} is chattering. The simulation of this solution can be seen in Figure \ref{nlinear}.

\section{Asymptotics}

Solutions of the  system  \eqref{mainc} admit   infinitely  many  jumps,  and   this makes, in general, impossible to  find an exact  solution  or  adequately to  simulate  it. So,  in this  section we suggest considering   degenerate equation to  find  the perturbed system  approximately.  In order  to  increase   the precision of  approximation we   follow the  idea of asymptotic   approximations.  Consider the system     
\begin{equation}\label{maincm}
\begin{split}
	\ddot{x}&=f(x,\dot{x}), \\
	\Delta \dot{x}&|_{\substack{x=\varphi \\ i < [ \frac{1}{r} ]}}=-(1+r)\dot{x}, 
\end{split}
\end{equation}
where $i$ is the index of impacts $\theta_i$, $[.]$ denotes the greatest integer function, with   additional   condition that   the number of impulsive moments   has to  be  not more than $\left[ \frac{1}{r} \right],$ i.e., $\theta_i, i=1,2,...,\left[ \frac{1}{r} \right].$ One can guarantee for the fixed value of the parameter $r$,  the incomplete chattering occurs only. The  number of impacts increases unboundedly as the parameter tends to zero. For this reason, we say that  system \eqref{maincm} admits the transient chattering. Assume that this system  satisfy all conditions of Theorem \ref{thmc}. For time $t>\theta_{\left[ \frac{1}{r} \right]},$ the system is only governed by $\ddot{x}=f(x,\dot{x}).$ Condition $(C1)$ implies that on the interval $[\theta_{\left[ \frac{1}{r} \right]},\theta_{\infty}]$ the bead stays on the position $x=\varphi.$ 
 
For each  its solution, system \eqref{maincm}  has  finite number  of discontinuity  moments.  That  is why,  one can  find an exact  solution of the problem or  at  least  it  is possible to  make proper simulations.    One can   easily see that  solutions of the last  system and system \eqref{mainc}   with   identical  initial data coincide   on the  interval  $[0,\theta_{\left[ \frac{1}{r} \right]}).$  They   are   different  only  in the interval $[\theta_{\left[ \frac{1}{r} \right]},\theta_{\infty}].$  The  length of the last  interval  diminishes to 0  as  $r\to 0.$ Consequently, the solutions of system \eqref{maincm} are asymptotic approximations for the solutions of system \eqref{mainc}.

\section{The Dynamics of Repeated Impacts against a Sinusoidally Vibrating Table}

In this section, we consider a mechanical model consisting of a bead bouncing  on a vibrating  table, which is investigated in the papers of Holmes and Guckenheimer \cite{holmes, guckenheimer}. It is demonstrated  that the model can generate chaos  \cite{holmes}. In this paper, we show that in the mechanism one can observe another type of complex dynamics, namely chattering.

Consider a bouncing bead colliding with a sinusoidally vibrating table. Assume that the table is so massive that it does not react to collisions with the bouncing bead and it moves according to law $X(t)=X_0\sin \omega t.$ The change of the velocity of the bouncing bead at the impact moment is given by the relation $r=\frac{\dot{X}_+-\dot{x}_+}{\dot{x}_--\dot{X}_-},$ where $r$ is the restitution coefficient, $0<r<1$, $\dot{X}_-, \dot{X}_+, \dot{x}_-, \dot{x}_+ $ are the velocities of the table and the bouncing bead before and after impact, respectively. Since the collision does not affect the velocity of the table, we can write $\dot{X}_-=\dot{X}_+.$ Then the model will be as follows
\begin{equation}
\begin{split}
\ddot{x}&=-g, \\
\Delta \dot{x}|_{x=X}&=-(1+r)(\dot{x}-\dot{X}),\\
X(t)&=X_0\sin(\omega t),
\end{split}
\label{beadontable}
\end{equation}
where $g$ is the gravitational acceleration ($g\approx  9.8 m/s^2$). 

Now, let us  consider a general form. Instead of gravitational constant $g,$  take a function $f(u,v).$ Then, the model will be of the form
\begin{equation}
\begin{split}
\ddot{x}&=f(x,\dot{x}), \\
\Delta \dot{x}|_{x=X}&=-(1+r)(\dot{x}-\dot{X}),\\
X(t)&=X_0\sin \omega t,
\end{split}
\label{table}
\end{equation}
where function $f(u,v)$ is a continuous function on the domain $G=\{X_0/10  \leq u \leq h, |v|\leq \bar{h} \},$ for fixed positive numbers $h, \bar{h},$ and it satisfies the local Lipschitz  condition in its variables on $G.$ Also, this system has conditions (C1), (C2) defined in the first section for all $(u,v)\in G$. By conditions on function $f(u,v)$ and compactness of the domain $G,$  we have a positive number $M$ such that $f(u,v)\geq -M$ for all $(u,v)\in G.$ 

\begin{figure}%
\includegraphics[width=\columnwidth]{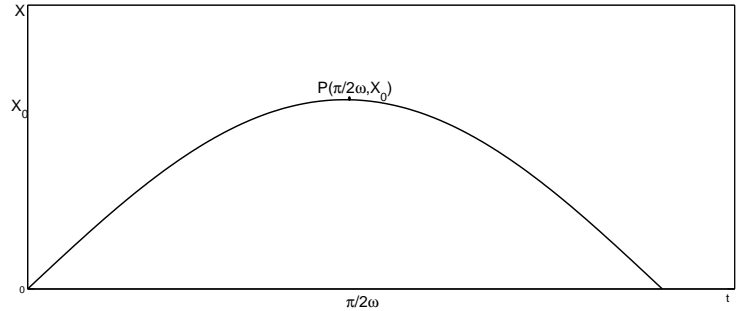}%
\caption{The graph of $X(t)=X_0\sin \omega t$ on the interval $[0,\pi/\omega].$}%
\label{vtable}%
\end{figure}

 Next, consider the graph of the function $X(t)=X_0\sin \omega t.$ The slope of the graph is $\dot{X}(t)=X_0\omega \cos \omega t.$ It is easily seen that if $\omega$ is small and $t$ is near $\pi/2\omega$ the graph is close to a horizontal line. Consequently, for sufficiently small $\omega$ and for time $t$ near $\pi/2\omega$ if the following inequality 
\begin{equation}\label{conditionhbar2}
M\sqrt{\frac{2(h-X_0/10)}{m}}<\bar{h},
\end{equation}
is true and conditions (C1), (C2) are satisfied, according to Theorem \ref{thmc}  there is  chattering for solutions whose integral curves are near to the point $P(\pi/2\omega,X_0)$ (see Figure \ref{vtable}). 

Finally, to demonstrate the result through simulation, we continue with the bouncing bead  on the sinusoidally vibrating table.

\begin{figure}[H]
\centering
\includegraphics[width=6cm]{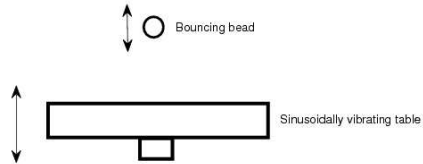}%
\caption{The bouncing bead on the sinusoidally vibrating table.}%
\label{balltable}%
\end{figure}

\begin{figure}[H]%
\includegraphics[width=12cm, height=6cm]{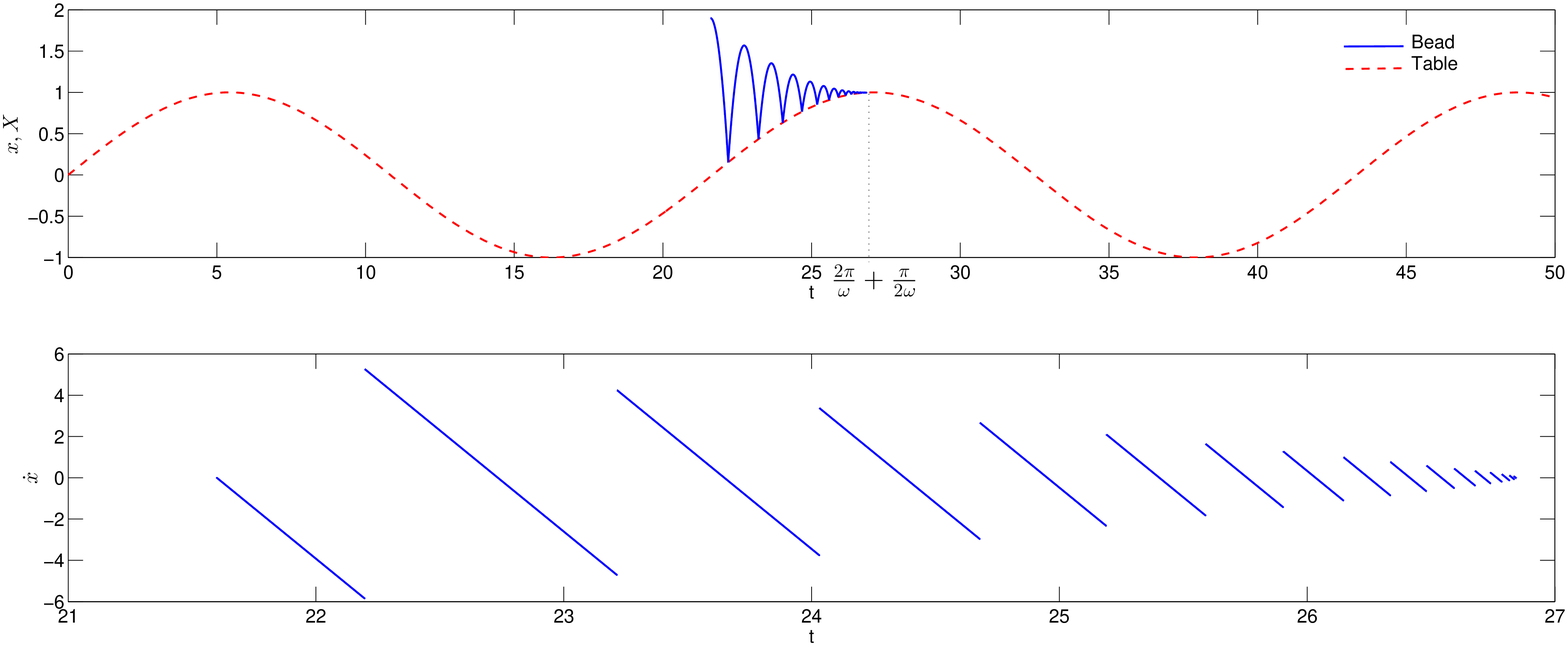}%
\caption{The graph of the coordinates of system \eqref{exc2}.}%
\label{bballtable2}%
\end{figure}

\textbf{Example 2. } Let us return to the bouncing bead  on the sinusoidally vibrating table, see Figure \ref{balltable}, with the same properties of system \eqref{table}. Then the model will be as follows
\begin{equation}
\begin{split}
\ddot{x}&=-g, \\
\Delta \dot{x}|_{x=X}&=-(1+r)(\dot{x}-\dot{X}),\\
X(t)&=X_0\sin(\omega t),
\end{split}
\label{exc2}
\end{equation}
where $t\geq 0.$ Let us take $\varphi=X_0=1$ and consider the domain $0.1\leq x\leq 2, |\dot{x}|<7.$ Then, we have  $f(x,\dot{x})=-g<0, $ $|f(x,\dot{x})|=|-g|=g=M=m$  and  $M\sqrt{\frac{2(h-X_0/10)}{m}}=\sqrt{37.24}<7.$ If we choose the initial conditions $x(2\pi/\omega)=1.9, \dot{x}(2\pi/\omega)=0$, where $r=0.9$, $\omega=0.29,$ it can be seen that the conditions of Theorem \ref{thmc} are satisfied and consequently, this solution is chattering. In Figure \ref{bballtable2}, one can observe the coordinates of system \eqref{exc2} which  supports our theoretical result.


\section{The Modified Moon-Holmes Model}
\begin{figure}[H]%
\centering
\includegraphics[width=8cm]{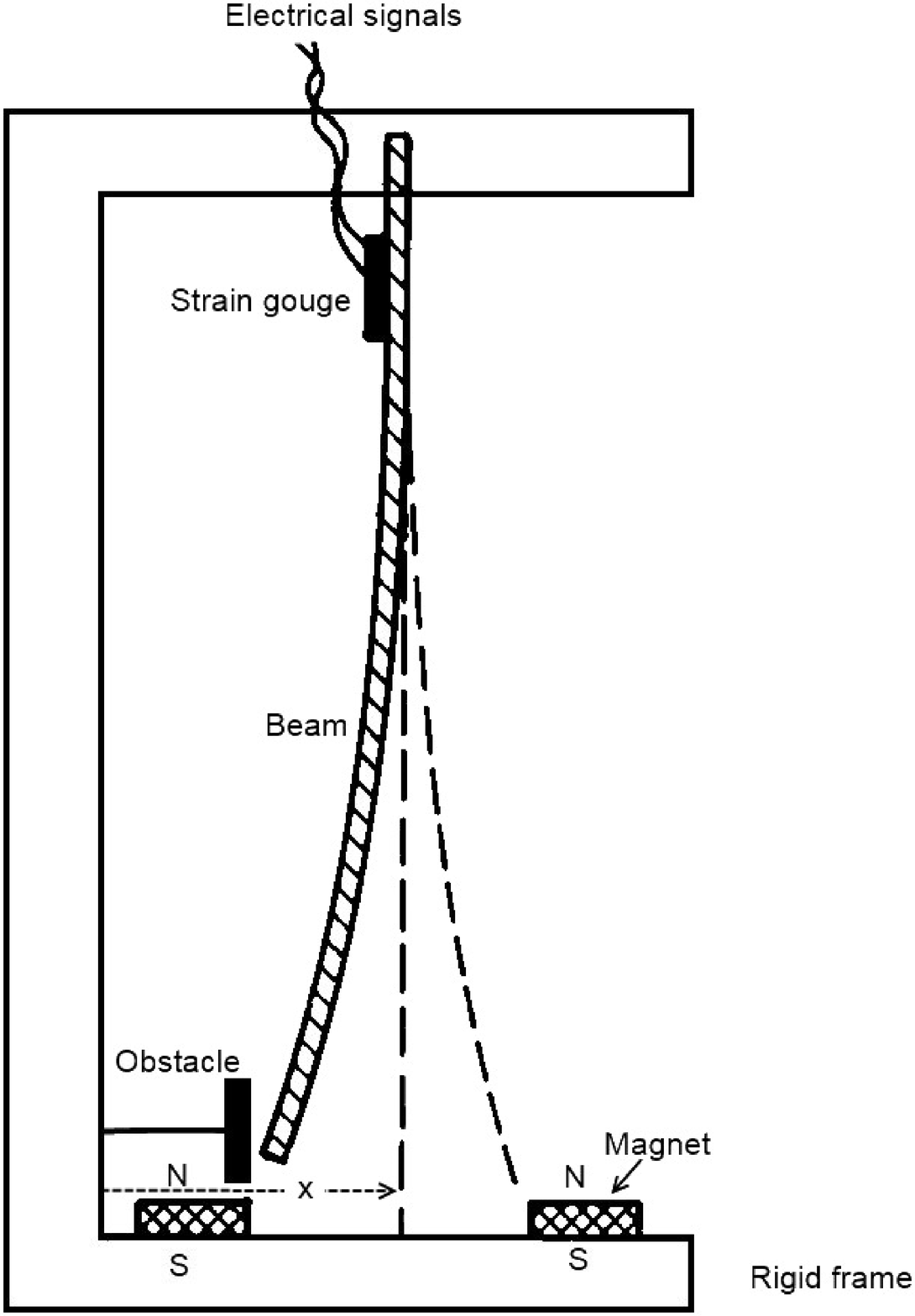}%
\caption{The magneto-elastic beam with the obstacle.}%
\label{duffing1}%
\end{figure}

The main task of this section is to consider the modified Moon-Holmes Model. Moon and Holmes \cite{moon1,moon2} showed that the Duffing equation in the form 
$$ \ddot{x} + \delta \dot{x} - x + x^3 = \gamma \cos wt  $$
provides the simplest possible model for the forced vibrations of a cantilever beam in the nonuniform field of two permanent magnets. Such an equation describes the dynamics of a buckled beam or plate when only one mode of vibration is considered. We modify the model as  adding a rigid obstacle over the magnet and in front of the beam such that the beam collides the obstacle  and from Newton Law of impacts it bounces back. (The system is sketched in Figure \ref{duffing1}.) The suggested model has the form of the following impulsive system
\begin{equation}
\begin{split}
\ddot{x}  &= - \delta \dot{x} + x - x^3+ \gamma \cos{wt}, \\
\Delta\dot{x}&|_{x=\varphi}=-(1+r)\dot{x},
\end{split}
\label{duffingimpulse1}
\end{equation}
where $x$ is the distance from the wall to the end of the beam, $\varphi$ is the position of the obstacle, $r$ is the restitution coefficient. Now, if  the coefficients $\gamma$ and $\delta$ are equal to zero, one obtains   

\begin{figure}[H]
\includegraphics[width=10cm, height=6cm]{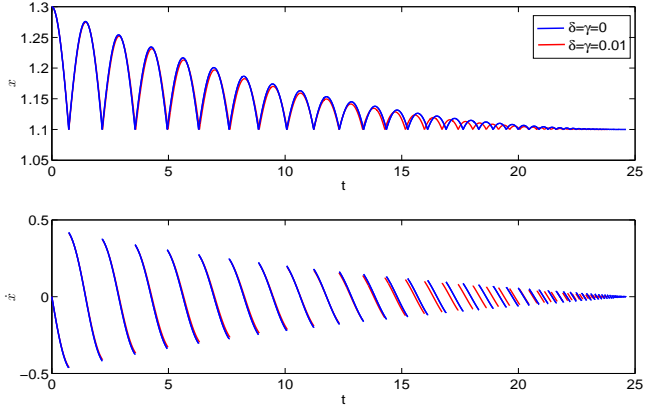}%
\caption{The coordinates of systems \eqref{duffingimpulse1} and \eqref{duffingimpulse2} with $w=0.1$. It can be seen that the solution of perturbed system \eqref{duffingimpulse1}  is also chattering.}%
\label{duffing2}%
\end{figure}

\begin{equation}
\begin{split}
\ddot{x} &= x - x^3, \\
\Delta\dot{x}&|_{x=\varphi}=-(1+r)\dot{x}.
\end{split}
\label{duffingimpulse2}
\end{equation}
For this system, choose $\varphi=1.1$ and for the domain $H$ let $h=1.5, \bar{h}=3.$ One can see that function $f(x,\dot{x})=x-x^3$ satisfies conditions (C1) and (C2), i.e. $-1.875\leq f(u,v)\leq -0.331$ for all $(u,v)\in H$ and $f(u,v)$ is an even function in $v.$ Moreover, condition \eqref{conditionhbar} is valid since $1.875\sqrt{\frac{2(1.5-1.1)}{0.331}}\approx 2,91<3.$ Therefore, by Theorem \ref{thmc} all solutions of system \eqref{duffingimpulse2} with initial values $(x(0),\dot{x}(0))=(x_0,0),$ $\varphi<x_0<h,$ are chattering. 
Obviously, system \eqref{duffingimpulse1} does not satisfy condition (C2). But, one can easily notice that for sufficiently small $\delta$ and $\gamma,$ by the continuous dependence on parameters and initial value for the impulsive differential equations with non-fixed moments \cite{Akhmet2010}, the solutions of \eqref{duffingimpulse1} with the same initial conditions of  \eqref{duffingimpulse2} are chattering as well. For the numerical simulation, let  $(x(0),\dot{x}(0))=(1.3,0)$ and $r=0.9$. Then, one can see that Figure \ref{duffing2} supports our theoretical discussion.


\section{Continuous Chattering}
In this section, we demonstrate the continuous chattering which is understood as infinitely many oscillations in finite time. Let us observe how continuous chattering appears if a mechanical model is perturbed with a discontinuous one. For this reason, we couple system \eqref{mainc} with the following equation of a   mass-spring-damper equation 
\begin{equation}\label{spring}
m\ddot{y}+c\dot{y}+ky=0,
\end{equation}
with mass $m$, spring constant $k$, and viscous damper of damping coefficient $c$. If the characteristic equation of system \eqref{spring} has roots with negative real parts, then it admits asymptotically stable equilibrium. By the argument of periodicity theorem for system with stable equilibrium, one can expect that in system \eqref{spring}, continuous chattering appears if it is perturbed by a chattering solution of \eqref{mainc}. Thus, let us write the coupled system taking $f(x,\dot{x})=-g$ in \eqref{mainc},  $m=1,c=3,k=2$ in \eqref{spring}  and $x=x_1,\dot{x}=x_2,y=x_3,\dot{y}=x_4,$ in the form
\begin{equation}\label{coupled}
\begin{split}
\dot{x}_1&=x_2,\\
\dot{x}_2&=-g,\\
\Delta &|_{x_1=1}=-(1+r)x_2 \\
\dot{x}_3&=x_4,\\
\dot{x}_4&=-2x_3-3x_4+20x_2^2.
\end{split}
\end{equation}
with initial conditions $x_1(0)=6, x_2(0)=0, x_3(0)=10, x_4(0)=-1000$ and $r=0.9.$
\begin{figure}[H] 
\includegraphics[width=\columnwidth]{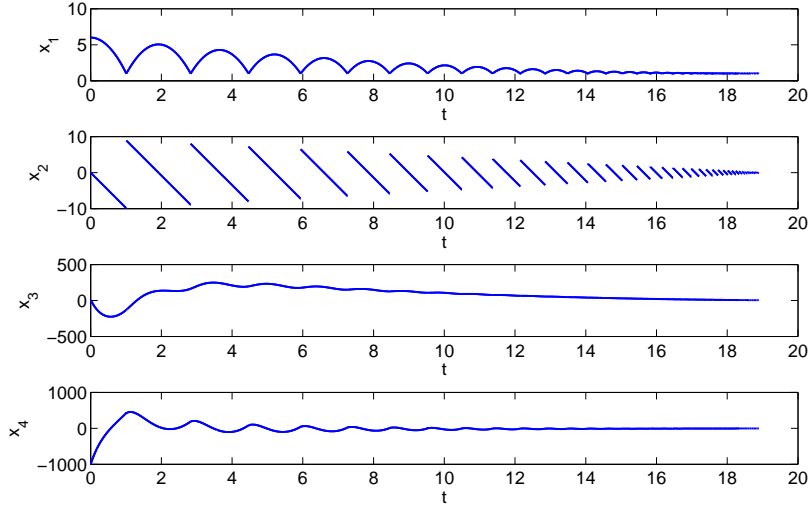}%
\caption{The graphs of the coordinates of  system \eqref{coupled}.}%
\label{a}%
\end{figure}
Since the first coupling is unilateral, the second equation does not influence the first one. That is why, its dynamics are the same as in Figure \ref{a}. But, for the second coupling in Figure \ref{a}, we can see the effect of perturbation which we call as continuous chattering.

\section{Pyragas Control}
There   are   many   papers  which are   searching methods to   minimize   and control   different  types of chattering \cite{pendulumcontrol, Chen,Alasty}.  The   problem  definitely    has to  be analyzed for   the impact  chattering, also. One can   accept   that  the control of   impact  chattering is  a   concrete   perturbation, which   brings the  system under control to  a regular motion.  That   is,   equilibria  or periodic motions.  In   the circumstances of the present research,  it  is desired that  a family   of chattering   solutions has to  be regularized, if not  all of them.   We   will  discuss, in this part of the paper,  system \eqref{exc}   of  Example $1.$ It  was  shown that any   solution of this system,  which  starts in a domain,  is chattering.   Let  us apply   the control of the form  $C[x_1(t-\tau)-x_1(t)]$  to   the system.   It  is applied, for instance,  to  stabilize periodic motions of chaotic dynamics, and it is called Pyragas control\cite{pyragas}. Now, we will apply the control   to  depress  the chattering   in the   system.  Let us construct the following  system denoting $x_1=x$ and $x_2=\dot{x}$ 
\begin{equation}\label{control1}
\begin{split}
\dot{x_1}&=x_2,\\
\dot{x_2}&=-x_1^3-\cos x_2+C[x_1(t-\tau)-x_1(t)] \\
\Delta x_2&|_{x_1=2}=-(1+r)x_2.
\end{split}
\end{equation}
 We performed a series of simulations of the system with fixed $C=-30,$ $\tau=1$ and $r=0.6.$ Consider $x_2(0)$ as it was requested to prove family of chattering solutions in Theorem \ref{thmc}. For the initial first coordinate $x_1(0)$ we tried values starting from  $2.5$ to $200$. 
\begin {table}[H]
\label{table1} 
\begin{center}
\begin{tabular}{l*{6}{|c}r}
$x_1(0)$ & 2.5 & 3 & 5 & 10 & 100 & 200  \\ \hline
Period $T$ & 1.22& 1.22& 1.22& 1.22& 1.22& 1.22 \\  \hline
Amplitude & 2.933 & 2.933& 2.933& 2.933& 2.933& 2.933 
\end{tabular}
\caption{Periods and amplitudes of the first coordinate $x_1(t)$ of system \eqref{control1} for different values of $x_1(0).$} 
\end{center}
\end{table}
For all these solutions the ultimate periodicity has been approved   with   period  $T=1.22.$ Observe that the period is different from the delay term $\tau=1.$ One can see from the table that the amplitudes are equal to $2.933$ for all values of $x_1(0)$ as well.  At the   same time,   the  chattering  has not  been   decaying   for the  solution with  $x_1(0)= 2.1.$  These all demonstrate   that    the control  problem  can   be solved for  the  chattering,  but certain  conditions    have   to  be  determined  to   specify   the controllable domains   and   conditions for the stability of the  arranged   periodic motions.  We  suppose   that   these problems will be researched  in  next papers.

\begin{figure}[H]
\centering
\includegraphics[height=5cm,width=10cm]{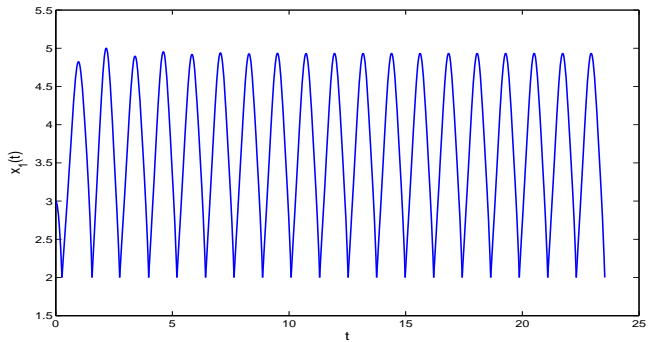}%
\caption{Simulation of the first coordinate $x_1(t)$ of controlled system \eqref{control1} with initial conditions $x_1(0)=3,x_2(0)=0$.}%
\label{controlfig1}%
\end{figure}

For  $x_1(0)=3,$ the  periodic orbit $x_1(t)$ can be seen in Figure \ref{controlfig1}, which shows the effectiveness of the control. 

\section{Conclusion}
In this article, we have considered the mechanical models  with impacts. For these models, the chattering phenomenon, which is defined    as a motion with  infinitely  many   discontinuities  in a finite time, is studied. The sufficient conditions are determined for the existence of the chattering. Asymptotics are discussed to find an approximation solution and to simulate the chattering solution. We study the famous example: the bouncing bead on a sinusoidally vibrating table which generates chaos \cite{holmes}. It is shown that this mechanism has chattering solutions. Furthermore, we modify the Moon-Holmes model \cite{moon1}, which yields chaos also, with an obstacle to obtain an impacting model. We demonstrate that this model provides chattering.
Perturbing a continuous mechanical process by a discontinuous one having chattering solutions, continuous chattering, which is defined as the appearance of infinitely many oscillations in a finite time, is constructed.

The application of  results  of paper \cite{nordmark} is to prove the existence of a unique chattering solution of the bouncing ball, see Section 3.5. At the same time, by simulation it is proven that a double pendulum admits chattering. Our method, in some sense, is more wider than the one in paper \cite{nordmark}. For example, by Theorem \ref{thmc} in this paper, we have verified that there are infinitely many chattering motions with initial values in an interval. Thus, the present result is complement to that one accomplished in \cite{nordmark}. However, our approach does not work for the double pendulum, since condition $(C2)$ is not valid for the model. Nevertheless, in our next investigation, we plan to extend the method without condition $(C2).$


\bibliography{ref}

\begin{thebibliography}{10}
\expandafter\ifx\csname url\endcsname\relax
  \def\url#1{\texttt{#1}}\fi
\expandafter\ifx\csname urlprefix\endcsname\relax\def\urlprefix{URL }\fi
\expandafter\ifx\csname href\endcsname\relax
  \def\href#1#2{#2} \def\path#1{#1}\fi

\bibitem{Lee}
H.~Lee, V.~I. Utkin, {Chattering suppression methods in sliding mode control
  systems}, Annual Reviews in Control 31~(2) (2007) 179--188.

\bibitem{ibrahim}
R.~A. Ibrahim, {Vibro-Impact Dynamics: Modeling, Mapping and Applications},
  Lecture Notes in Applied and Computational Mechanics, Springer, 2009.

\bibitem{Giusepponi}
S.~Giusepponi, F.~Marchesoni, M.~Borromeo, {Randomness in the bouncing ball
  dynamics}, Physica A: Statistical Mechanics and its Applications 351~(1)
  (2005) 142--158.

\bibitem{budd}
C.~Budd, F.~Dux, {Chattering and Related Behaviour in Impact Oscillators},
  Philosophical Transactions: Physical Sciences and Engineering 347~(1683)
  (1994) pp. 365--389.

\bibitem{nordmark}
A.~B. Nordmark, P.~T. Piiroinen, {Simulation and stability analysis of
  impacting systems with complete chattering}, Nonlinear Dynamics 58~(1-2)
  (2009) 85--106.

\bibitem{Arnold}
R.~N. Arnold, {Cutting Tools Research: Report of Subcommittee on Carbide Tools:
  The Mechanism of Tool Vibration in the Cutting of Steel}, Proceedings of the
  Institution of Mechanical Engineers 154~(1) (1946) 261--284.

\bibitem{nagaev}
R.~F. Nagaev, E.~B. Kremer, {Mechanical Processes with Repeated Attenuated
  Impacts}, World Scientific, 1999.

\bibitem{Kob}
A.~A. Kobrinsky, A.~E. Kobrinsky, {Two Dimensional Vibro-Impact
  Systems(Russian)}, Nauka, 1981.

\bibitem{akhmetayse}
M.~U. Akhmet, A.~Kivilcim, {The Models with Impact Deformations},
  Discontinuity, Nonlinearity, and Complexity 4~(1) (2015) 49--78.

\bibitem{holmes}
P.~J. Holmes, {The dynamics of repeated impacts with a sinusoidally vibrating
  table}, Journal of Sound and Vibration 84~(2) (1982) 173--189.

\bibitem{guckenheimer}
J.~Guckenheimer, P.~Holmes, {Nonlinear oscillations, dynamical systems, and
  bifurcations of vector fields}, Applied mathematical sciences,
  Springer-Verlag, 1990.

\bibitem{moon1}
F.~C. Moon, P.~J. Holmes, {A magnetoelastic strange attractor}, J. Sound Vib.
  65~(2) (1979) 285--296.

\bibitem{moon2}
F.~C. Moon, P.~J. Holmes, {Addendum: a magnetoelastic strange attractor}, J.
  Sound Vib. 69~(2) (1980) 339.

\bibitem{Akhmet2010}
M.~Akhmet, {Principles of Discontinuous Dynamical Systems}, Springer, New York,
  2010.

\bibitem{pendulumcontrol}
I.~Boussaada, I.-C. Morarescu, S.-I. Niculescu, {Inverted Pendulum
  Stabilization Via a Pyragas-Type Controller: Revisiting the Triple Zero
  Singularity}, in: Preprints of the 19th World Congress The International
  Federation of Automatic Control, 2014, pp. 6806--6811.

\bibitem{Chen}
W.-C. Chen, {Dynamics and control of a financial system with time-delayed
  feedbacks}, Chaos, Solitons {\&} Fractals 37~(4) (2008) 1198--1207.

\bibitem{Alasty}
A.~Alasty, H.~Salarieh, {Controlling the chaos using fuzzy estimation of
  {\{}OGY{\}} and Pyragas controllers}, Chaos, Solitons {\&} Fractals 26~(2)
  (2005) 379--392.

\bibitem{pyragas}
K.~Pyragas, {Continuous control of chaos by self-controlling feedback}, Physics
  Letters A 170~(6) (1992) 421--428.

\end{thebibliography}

\end{document}